\documentclass{amsart}

\newcommand{\remark}{\noindent {\bf Remark. }}
\newcommand{\corollary}{\noindent {\bf Corollary. }}
\newcommand{\ws}{\hspace{4pt}}
\newtheorem{theorem}{Theorem}

\newtheorem{proposition}{Proposition}
\newtheorem{lemma}{Lemma}
\newtheorem{defi}{Definition}
\usepackage[]{amssymb,amsopn}
\usepackage[]{amsmath}
\usepackage[]{eucal}
\usepackage[]{latexsym}

\begin{document}

\title[Energy function]{The energy function with respect to the zeros of the exceptional Hermite polynomials}
\author{\'A. P. Horv\'ath }

\thanks{Supported by Hungarian National Foundation for Scientific Research, Grant No.  K-100461}

\begin{abstract}
We examine the energy function with respect to the zeros of exceptional Hermite polynomials. The localization of the eigenvalues of the Hessian is given in the general case. In some special arrangements we have a more precise result on the behavior of the energy function. Finally we investigate the energy function with respect to the regular zeros of the exceptional Hermite polynomials.

\end{abstract}
\maketitle

\section{Introduction}

Exceptional orthogonal polynomials were introduced recently by D. G\'omez-Ullate, N. Kamran and R. Milson (cf. e.g. \cite{gukm3} and \cite{gukm2}). Besides the classical ones, exceptional orthogonal polynomials play a fundamental role in the construction of bound-state solutions to exactly solvable potentials in quantum mechanics. The relationship between exceptional orthogonal polynomials and the Darboux transform is observed by C. Quesne (cf. e.g. \cite{q}). Higher-codimensional families were introduced by S. Odake and R. Sasaki \cite{os}. The last few years have seen a great deal of activity in the area of exceptional orthogonal polynomials (cf. e.g. \cite{d},\cite{d1}, \cite{ggm}, \cite{gukkm} and the references therein).
Exceptional orthogonal polynomials are complete orthogonal systems with respect to a weight on an interval $I$, and they are also eigenfunctions of a second order differential operator. In all known cases the above mentioned weight is
$$ w(x)=\frac{w_0(x)}{P^2(x)},$$
where the classical counterpart of the exceptional polynomials in question are orthogonal with respect to $w_0$ on $I$. $P(x)$ is a polynomial with zeros in the exterior of  $I$. The sequence of the degrees of exceptional orthogonal polynomials are gappy, and the number of the missing degrees depends on the degree of $P$.
The zeros of exceptional orthogonal polynomials are divided into two groups: regular zeros which lie in the domain of orthogonality, and exceptional zeros wich lie in the exterior of $I$. In connection with the location of zeros of exceptional orthogonal polynomials the following conjecture is drafted in \cite{km}:

"The regular zeros of exceptional orthogonal polynomials have the same asymptotic behavior as the zeros of their classical counterpart. The exceptional zeros converge to the zeros of $P(x)$."

The location of zeros of exceptional Laguerre and Jacobi polynomials are described by D. G\'omez-Ullate, F. Marcell\'an and R. Milson \cite{gumm}, of exceptional Hermite polynomials by A. Kuijlaars and R. Milson, \cite{km}. Now we can draw up the following conjecture:

The normalized counting measure on the scaled regular zeros of exceptional orthogonal polynomials tend in the weak star topology to the equilibrium measure defined by $w_0$.

We investigate here the energy function at the zeros of exceptional Hermite polynomials.

At the end of the nineteenth century the maximum of a function of $n$ variables like
$T(u_1,\dots , u_n)=\prod_{1\leq i < j\leq n}(u_i-u_j)^2$
was already investigated by some authors e.g. Stieltjes and Hilbert.
Few years later I. Schur proved that on finite or infinite intervals (under certain conditions) the maximum is attained at the zeros of certain classical orthogonal polynomials (cf.  \cite{sch} and the references therein). The question also has an electrostatic interpretation. If a "nice" weight function ($w$) is responsible for the external field, the energy function above will change to
$$T_w(u_1,\dots , u_n)=\prod_{i=1}^nw(u_i)\prod_{1\leq i < j\leq n}(u_i-u_j)^2.$$
According to the usual notations introduced in \cite{st}
$$\inf_{u_1,\dots , u_n\in H}-\frac{1}{n(n-1)}\log T_w(u_1,\dots , u_n)$$
is the discrete energy or the $n^{th}$ transfinite diameter of the set $H$ with respect to $w^{\frac{1}{2(n-1)}}$. These systems of minimal energy are expansively examined and have several applications, for instance these are the optimal systems of nodes of interpolation, see e.g. \cite{f}, \cite{aho}-\cite{i}. If $w$ is a classical weight function, then the solutions of the weighted energy problem are the sets of the zeros of classical orthogonal polynomials, namely with Jacobi weights $w_{\alpha,\beta}=(1-x)^{\alpha}(1+x)^{\beta}$, with Laguerre weights $w_{\alpha}=x^{\alpha}e^{-x}$, with Hermite weight $w=e^{-x^2}$, the zeros of orthogonal polynomials with respect to $w_{\alpha-1,\beta-1}$,  $w_{\alpha-1}$, $w$ respectively (cf. e. g. \cite{aho} and the references therein).

One can ask whether the behavior of the zeros of the exceptional polynomials are similar to the classical ones.
The electrostatic interpretation of zeros of $X_1$-Jacobi polynomials is given by D. Dimitrov and Yen Chi Lun \cite{dy}, and of $X_m$-Jacobi and $X_m^I$-Laguerre, $X_m^{II}$-Laguerre polynomials see \cite{ahor}. In case of exceptional Hermite polynomials the situation is rather different, since the exceptional zeros must be complex. We formulate the question as a real maximalization problem. (The other possibilities are drawn up in the remarks in the next section.) We find that if we introduce $2m+n$ real variables (the real and the imaginary parts of the $m$ exceptional zeros and the $n$ real zeros), then we will get a result similar to the exceptional Laguerre and Jacobi cases, and if we introduce $m+n$ real variables (the real parts of the $m$ exceptional zeros and the $n$ real zeros), then we will get a result similar to the classical Hermite case. Finally we show that the set of the regular zeros is the optimal system with respect to certain varying weights.

\section{Definitions and Results}

\subsection{Exceptional Hermite polynomials}
Following the notations of \cite{km} a partition $\lambda$ of length $r$ is $$\lambda=(\lambda_1\geq \lambda_2\geq \dots \geq \lambda_r\geq 1).$$  The size of $\lambda$ is $|\lambda|=\sum_{i=1}^r\lambda_i$. The generalized Hermite polynomial associated with $\lambda$ is
 $$H_{\lambda}=\mathrm{Wr}\left[H_{\lambda_r},\dots,H_{\lambda_2+r-2},H_{\lambda_1+r-1}\right],$$
where  $\mathrm{Wr}$ denotes the Wronskian determinant. The degree of $H_{\lambda}$ is
$$\deg H_{\lambda}=|\lambda|.$$
We call $\lambda$ an even partition if $r$ is even and $\lambda_{2k-1}=\lambda_{2k}$ for $k=1, \dots , \frac{r}{2}$. In this case $H_{\lambda}$ has no zeros on the real axis cf. \cite{a}, \cite{ggm} and \cite{k}. Let
$$w_{\lambda}=\frac{e^{-x^2}}{\left(H_{\lambda}(x)\right)^2}$$
a positive weight function on the real line. For a fixed partition $\lambda$ let
$$N_{\lambda}:=\{n \geq |\lambda|-r : n\neq |\lambda|+\lambda_i-i, \ws i=1,\dots , r\},$$
and for $n \in N_{\lambda}$ let
$$P_n:=\mathrm{Wr}\left[H_{\lambda_r},\dots,H_{\lambda_2+r-2},H_{\lambda_1+r-1}, H_{n-|\lambda|+r}\right].$$
Then $\deg P_n=n$ and $P_n$ satisfies the differential equation (cf. \cite{ggm})
\begin{equation}P_n''+\left(-2x-2\frac{H_{\lambda}'}{H_{\lambda}}\right)P_n'+\left(\frac{H_{\lambda}''}{H_{\lambda}}+2x\frac{H_{\lambda}'}{H_{\lambda}}+2n-|\lambda|\right)P_n=0,\end{equation}
that is $\frac{P_n(x)e^{-\frac{x^2}{2}}}{H_{\lambda}(x)}$ is an eigenfunction of the differential operator
\begin{equation} -y''+\left(x^2-2\left(\log H_{\lambda}\right)''\right)y,\end{equation}
with eigenvalue $2n-2|\lambda|+1$. Moreover if $\lambda$ is an even partition and $n,m \in N_{\lambda}$, $n\neq m$
$$\int_{\mathbb{R}}P_nP_mw_{\lambda}=0.$$
The closure of the span of $\{P_n\}_{n \in N_{\lambda}}$ is dense in the Hilbert space $L^2(\mathbb{R},w_{\lambda})$, cf. \cite{d}.

Throughout the rest of the paper, $\lambda$ is a fixed even partition. For simplicity let $ H_{\lambda}=H$, $\deg H =m$. We assume that all the zeros of $H$ are simple that is $H(z)=c\prod_{k=1}^m(z-w_k)$, where $w_k=u_k+iv_k$. The condition of simple zeros may not be too restrictive since by a conjecture (cf. \cite{fhv}) for any partition $\lambda$ the zeros of $H_{\lambda}$ are simple, except possibly for the zero at $z=0$.
The conjecture is known  in the special case $\lambda=(\nu,\nu)$, cf. \cite[Lemmas 3.1 and 3.2]{gfgu}.

\medskip

\subsection{Energy function}
Analogously to the Laguerre and Jacobi cases, let $w(z)=\frac{\left|e^{-z^2}\right|}{|H(z)|^2}$ on $\mathbb{C}$. This function is responsible for the "external field". It is positive, but is not bounded around $w_k$-s. It has "saddle points" and the behavior at the boundary of the domain is not so nice as in the classical cases, even if the domain is restricted to a strip $|\Im z|\leq c(m)$. The properties on the real line are much better. Indeed if $x\in \mathbb{R}$, $w(x)=\frac{e^{-x^2}}{H^2(x)}$, which is positive and bounded on the real line, since $H$ has no zeros on $\mathbb{R}$ and it tends to zero when $|x|$ tends to infinity.
We investigate the energy function at the zeros of exceptional Hermite polynomials with respect to $w$. Since the behavior of the energy function is interesting in case of large $n$, subsequently we assume that $n\geq \lambda_1$. As in \cite{ahor} we denote the exceptional Hermite polynomials by $P_{m,m+n}$, where $m+n\in N_{\lambda}$, and it shows that it has $m$ exceptional and $n$ regular, real zeros cf. \cite[(2.9)]{km}. By \cite[Theorem 2.3]{km} the exceptional zeros are not real (if $n$ is large enough). Later (if $n$ is large enough) we use the decomposition $P_{m,m+n}=P_mq_n$, where the zeros of $P_m$, $\{z_k=\xi_k+i\eta_k, k=1, \dots , m\}$ are the exceptional zeros of $P_{m,m+n}$, and the zeros of $q_n$, $\{x_i, i=1, \dots, n \}$ are the regular (real) zeros of $P_{m,m+n}$. Let us denote by $Z$ the set of the zeros of $P_{m,m+n}$.

The energy function is as follows. Let $U=u_1,\dots ,u_{m+n}\subset \mathbb{C}$.
$$T_w(U)=\prod_{k=1}^{m+n}w(u_k)\prod_{1\leq k<l\leq m+n}|u_k-u_l|^2.$$
For simplicity let us denote the logarithm of the energy function by $F$. We investigate $F$ as a function of $2m+n$ real variables.
$$F(\dot{\xi}_1,\dot{\eta}_1, \dots ,\dot{\xi}_m,\dot{\eta}_m,\dot{x}_1,\dots ,\dot{x}_n)=\log T_w=\sum_{k=1}^m\log w(\dot{z}_k)+\sum_{i=1}^n\log w(\dot{x}_i)$$ $$+\sum_{1\leq k<l\leq m}\log |\dot{z}_k-\dot{z}_l|^2+\sum_{1\leq i<j\leq n}\log |\dot{x}_i-\dot{x}_j|^2+\sum_{k=1}^m\sum_{i=1}^n \log |\dot{x}_i-\dot{z}_k|^2.$$

\medskip

\remark
Obviously 
$$F=\Re F_c,$$
where
$$F_c(u_1, \dots,u_{m+n})=-\sum_{i=1}^{m+n}u_i^2-2\sum_{i=1}^{m+n}\log H(u_i)+2\sum_{1\leq i<j\leq m+n}\log(u_i-u_j).$$
Recalling, that $H=H_{\lambda}$, by the differential equation (1) it is clear that all the first partial derivatives of $F$ are zero at $Z$.

\medskip

Let us denote the Hessian of $F$ by ${\bf H}\in \mathbb{R}^{2m+n\times 2m+n}$. Now we can state

\medskip

\begin{theorem}

If $n$ is large enough, ${\bf H}(Z)$ is nonsingular.\end{theorem}

\medskip

Actually we prove that ${\bf \hat{H}}=D^{-1}{\bf H}D$ is nonsingular, where $D$ is a diagonal matrix such that $d_{kk}=1$, $k=1, \dots, 2m$, $d_{2m+i2m+i}=K(m)$, $i=1,\dots n$. Here $K(m)$ is a constant depends on $m$ and is given later.

\medskip

Next we investigate the Gersgorin set $G$ of the partitioned matrix ${\bf \hat{H}}$. According to \cite[Theorem 6.3]{v} the eigenvalues of  ${\bf \hat{H}}$ are in $G$. $G$ consists of two subsets, $G=G_r\cup G_e$, where $G_r$ is generated by the regular zeros (see (6)), and $G_e$ is generated by the exceptional zeros (see (7)). Since ${\bf H}$ is symmetric, the eigenvalues of ${\bf \hat{H}}$ are real, thus we mean $G_r=G_r\cap\mathbb{R}$, and $G_e=G_e\cap\mathbb{R}$. (For the precise definition of the Gersgorin set of a partitioned matrix see the next section.) Then we have

\medskip

\begin{theorem}
Let $n$ be large enough. $G_r\subset \mathbb{R}_-$, and if $x\in G_r$ then $|x|\leq c n$. If $x\in G_e$, then $|x|\sim n$.\end{theorem}

\medskip

\remark
The other possibility is to investigate the energy function as a function of $m+n$ complex variables. In this case the Hessian would be ${\bf H}_{ij}=\frac{\partial^2 F}{\partial u_i\partial \bar{u}_j}$, and so ${\bf H}_{ij}(Z)=0$, $1\leq i,j \leq m$. This shows that $w=\frac{|e^{-z^2}|}{(\cdot)}$ has saddle points thus it is not a "nice" weight function. The usual extension of $e^{-x^2}$ as a weight to the complex plane is  $e^{-|z^2|}$. In this case the Hessian would be nice, ${\bf H}_{ij}(Z)=0$, $1\leq i,j \leq m, \ws i\neq j$ and ${\bf H}_{ii}(Z)=-1$, $1\leq i \leq m$, but the first partials are not zero at $Z$.

\medskip

Let $\{h_n\}_{n=0}^{\infty}$ be the sequence of orthonormal Hermite polynomials. Now we examine the special case when
$$H=d_{\nu}:= \left|\begin{array}{cc}h_{\nu}&h_{\nu+1}\\h'_{\nu}&h'_{\nu+1}\end{array}\right|.$$
By Christoffel-Darboux formula (cf. e.g. \cite[3.2.4]{sz}) $d_{\nu}=\frac{\gamma_{n+1}}{\gamma_n}\sum_{k=0}^{\nu}h_k^2$, ($\gamma_l>0$) is obviously positive on the real axis and let us recall that by \cite[Proposition 3.2]{gfgu} it has simple roots. Moreover the location of zeros is symmetric since $d_{\nu}$ is an even polynomial of degree $2\nu=m$ and it has real coefficients. In this special case the behavior of the energy function with respect to the zeros of $P_{m,m+n}$ is similar to the Jacobi and Laguerre cases (cf. \cite[Theorem 2]{dy}, \cite[Theorem 1]{ahor}).

\medskip

\begin{theorem}
Let $H=d_{\nu}$. If $n$ is large enough, the logarithmic energy function $F=\log T_w$ has a "saddle point" at the set of zeros of orthogonal polynomials with respect to $w$, more precisely the modified Hessian ${\bf \hat{H}}$ is diagonally dominant, $h_{2m+i2m+i}<0$, $i=1,\dots ,n$ and $h_{2l-1,2l-1}=-h_{2l,2l}< 0$, $l=1, \dots , m$.\end{theorem}

\medskip

\remark
The Hermite weight is exceptional in point of energy function, since usually the zeros of the orthogonal polynomials with respect to a weight $w$, are not the optimal systems with respect to the minimal energy problem with the external field generated by $w$ itself, except when $w$ is the Hermite weight (cf. e.g. \cite{sz}, \cite{aho}, \cite{i}, \cite{ho}). As it is shown in that special case above, on one hand the energy function at the zeros of exceptional Hermite polynomial behaves like in exceptional Jacobi and Laguerre cases, but on the other hand if we fixed $m$ variables at the imaginary parts of the exceptional zeros, the new function with $m+n$ variables has a maximum at $\xi_1,\dots ,\xi_m,x_1,\dots , x_n$ which seems to preserve the above mentioned property of the classical Hermite weight.

\medskip

Since the zeros of the exceptional Hermite polynomials are not the optimal set with respect to the energy problem generated by $w$, as in \cite{ho} and \cite{aho} we give a new weight $w_1$ on $\mathbb{R}$ which depends on $n$ and for which the regular zeros of the exceptional Hermite polynomials are optimal.

\medskip

\begin{defi}
A nonnegative weight is approximating on $(a,b)\subset \mathbb{R}$, if it has finite moments, it is twice differentiable and $\left(\log\left(\frac{1}{w}\right)\right){''}\geq 0$ on $(a,b)$, and if $a$ is finite, then $\lim_{ x\to a+ }\frac{w(x)}{x-a}=0$, and if $b$ is finite, then $\lim_{ x\to b-} \frac{w(x)}{b-x}=0$.
\end{defi}

\medskip

Let us take into consideration the general form of the differential equation of general orthogonal polynomials cf. \cite{m}, \cite{iw}. Reformulating and summarizing the results in \cite{ahor} we have the following

\begin{proposition} Let $q(x)=\prod_{k=1}^n(x-x_k)$ a polynomial with different real zeros which fulfils the following differential equation
$$q''+M_nq'+N_nq=0.$$
If
$$w=w_n=e^{\int M_n}$$
is an approximating weight, then $Z=\{x_1, \dots , x_n\}$ is the unique solution of the energy problem with respect to $T_w$.
\end{proposition}

\medskip

\begin{proposition}
Let the ordinary (orthogonal) polynomials of degree $n$ fulfil the following differential equation
$$y''+M_{0,n}y'+N_{0,n}y=0,$$
and let us denote by $w_0=w_{0,n}=e^{\int M_{0,n}}$.
Let us assume that the corresponding exceptional polynomials fulfil the modified differential equation
$$y''+M_ny'+N_ny=0,$$
where
$$M_n=M_{0,n}-2\frac{H'}{H},$$
where $H$ is a polynomial (of degree m) with simple zeros, cf.\cite[(1.1)]{km}. Let $P_{m,m+n}=P_mq_n$ as above, and let $w_1=w_{1,n}=e^{\int M_{1,n}}$, where $M_{1,n}=M_n+2\frac{P_m'}{P_m}$. If the exceptional zeros of the exceptional polynomials tend to the zeros of $H$ when $n$ tends to infinity, then $w_{1,n}$ is approximating if $n$ is large enough, and $w_0$ is approximating.
In this case the logarithmic energy function with respect to the weight $w_1$ attains its maximum at a unique set, which is the set of the regular zeros of $P_{m,m+n}$.
\end{proposition}

\medskip

\corollary
If $n$ is large enough, the logarithmic energy function with respect to the weight $w_1:= wP_m^2$ attains its maximum at a unique set, which is the set of the regular zeros $x_1, \dots , x_n$ of the exceptional Hermite polynomial $P_{m,m+n}$.

\medskip

\section{Proofs}

First we compute the elements of the Hessian. These element also can be formulated as the real or imaginary parts of the second partials of $F_c$, but as it is remarked it does not useful for our purposes, so we give the elements in real form.
$$h_{2k-1,2k-1}=F''_{\dot{\xi}_k\dot{\xi}_k}(Z)=-2 + 2\sum_{l=1}^m\frac{(\xi_k-u_l)^2-(\eta_k-v_l)^2}{((\xi_k-u_l)^2+(\eta_k-v_l)^2)^2}$$ $$-2\sum_{1\leq l\leq m \atop l\neq k}\frac{(\xi_k-\xi_l)^2-(\eta_k-\eta_l)^2}{((\xi_k-\xi_l)^2+(\eta_k-\eta_l)^2)^2}-2\sum_{i=1}^n\frac{(\xi_k-x_i)^2-\eta_k^2}{((\xi_k-x_i)^2+\eta_k^2)^2},$$
$$h_{2k-1,2k}=F''_{\dot{\eta}_k\dot{\xi}_k}(Z)= 4\sum_{l=1}^m\frac{(\xi_k-u_l)(\eta_k-v_l)}{((\xi_k-u_l)^2+(\eta_k-v_l)^2)^2}$$ $$-4\sum_{1\leq l\leq m \atop l\neq k}\frac{(\xi_k-\xi_l)(\eta_k-\eta_l)}{((\xi_k-\xi_l)^2+(\eta_k-\eta_l)^2)^2}-4\sum_{i=1}^n\frac{(\xi_k-x_i)\eta_k}{((\xi_k-x_i)^2+\eta_k^2)^2},$$
$$h_{2k-1,2l-1}=F''_{\dot{\xi}_l\dot{\xi}_k}(Z)=2\frac{(\xi_k-\xi_l)^2-(\eta_k-\eta_l)^2}{((\xi_k-\xi_l)^2+(\eta_k-\eta_l)^2)^2},$$
$$h_{2k-1,2l}=F''_{\dot{\eta}_l\dot{\xi}_k}(Z)=4\frac{(\xi_k-\xi_l)(\eta_k-\eta_l)}{((\xi_k-\xi_l)^2+(\eta_k-\eta_l)^2)^2},$$
$$h_{2k-1,2m+i}=F''_{\dot{x}_i\dot{\xi}_k}(Z)=2\frac{(\xi_k-x_i)^2-\eta_k^2}{((\xi_k-x_i)^2+\eta_k^2)^2}.$$
and
$$h_{2k,2k}=F''_{\dot{\eta}_k\dot{\eta}_k}(Z)=-F''_{\dot{\xi}_k\dot{\xi}_k}(Z)=-h_{2k-1,2k-1},$$ $$h_{2k,2k-1}=F''_{\dot{\xi}_k\dot{\eta}_k}(Z)=F''_{\dot{\eta}_k\dot{\xi}_k}(Z)=h_{2k-1,2k},$$

$$h_{2k,2l}=F''_{\dot{\eta}_l\dot{\eta}_k}(Z)=-F''_{\dot{\xi}_l\dot{\xi}_k}(Z)=-h_{2k-1,2l-1},$$ $$  h_{2k,2l-1}=F''_{\dot{\xi}_l\dot{\eta}_k}(Z)=F''_{\dot{\eta}_l\dot{\xi}_k}(Z)=h_{2k-1,2l},$$

$$h_{2k,2m+i}=F''_{\dot{x}_i\dot{\eta}_k}(Z)=4\frac{(\xi_k-x_i)\eta_k}{((\xi_k-x_i)^2+\eta_k^2)^2}.$$

$$h_{2m+i,2m+i}=F''_{\dot{x}_i\dot{x}_i}(Z)=-2 +2\sum_{k=1}^m\frac{(x_i-u_k)^2-v_k^2}{((x_i-u_k)^2+v_k^2)^2}$$ $$-2\sum_{k=1}^m\frac{(\xi_k-x_i)^2-\eta_k^2}{((\xi_k-x_i)^2+\eta_k^2)^2}-2\sum_{1\leq j \leq n \atop j\neq i}\frac{1}{(x_i-x_j)^2},$$
where the imaginary part of $\left(\frac{H'}{H}\right)'(x_i)$ is disappeared, because $H$ has real coefficients thus it has conjugate pairs of roots.
$$h_{2m+i,2k-1}=F''_{\dot{\xi}_k\dot{x}_i}(Z)=F''_{\dot{x}_i\dot{\xi}_k}(Z)=h_{2k-1,2m+i},$$
$$h_{2m+i,2k}=F''_{\dot{\eta}_k\dot{x}_i}(Z)=F''_{\dot{x}_i\dot{\eta}_k}(Z)=h_{2k,2m+i},$$
$$h_{2m+i,2m+j}=\frac{2}{(x_i-x_j)^2}.$$

For the proof of Theorem 1 we need some lemmas. In the next lemma, and the rest of the paper the matrix norm is $l_{\infty}$-norm, that is if $A \in \mathbb{C}^{n\times m}$, then $\|A\|=\max_{i=1}^n\left(\sum_{j=1}^m|a_{i,j}|\right)$.

Let $$A_{k,k}=\left[\begin{array}{cc}F''_{\dot{\xi}_k\dot{\xi}_k}(Z)&F''_{\dot{\eta}_k\dot{\xi}_k}(Z)\\F''_{\dot{\xi}_k\dot{\eta}_k}(Z)&F''_{\dot{\eta}_k\dot{\eta}_k}(Z)\end{array}\right].$$

\medskip

\lemma

There is a positive constant $c$ such that
\begin{equation}\|A_{k,k}^{-1}\|^{-1}> c n.\end{equation}

$$\sum_{1\leq l \leq m \atop l\neq k}\max\left\{\left|h_{2k-1,2l-1}\right|+\left|h_{2k-1,2l}\right|;\left|h_{2k,2l-1}\right|
+\left|h_{2k,2l}\right|\right\}$$
\begin{equation}+\sum_{i=1}^n\max\left\{\left|h_{2k-12m+i}\right|; \left|h_{2k2m+i}\right|\right\}=O(\sqrt{n}).\end{equation}

\medskip

\proof

To prove (3) first let us observe that the following two sums can be estimated by a constant which depends only on $m$:
$$S_1:=\sum_{1\leq l\leq m \atop l\neq k}\left(\frac{(\xi_k-u_l)^2-(\eta_k-v_l)^2}{((\xi_k-u_l)^2+(\eta_k-v_l)^2)^2}-\frac{(\xi_k-\xi_l)^2-(\eta_k-\eta_l)^2}{((\xi_k-\xi_l)^2+(\eta_k-\eta_l)^2)^2}\right)\leq c(m),$$ $$S_2:=2\sum_{1\leq l\leq m \atop l\neq k}\left(\frac{(\xi_k-u_l)(\eta_k-v_l)}{((\xi_k-u_l)^2+(\eta_k-v_l)^2)^2}-\frac{(\xi_k-\xi_l)(\eta_k-\eta_l)}{((\xi_k-\xi_l)^2+(\eta_k-\eta_l)^2)^2}\right) \leq c(m).$$
(Actually, since by \cite[Theorem 2.3.]{km} the exceptional zeros of $P_{m,m+n}$ tends to the zeros of $H$, the two sums above are small when $n$ is large.)
$$|\det A_{k,k}|=\left|-(F''_{\dot{\xi}_k\dot{\xi}_k}(Z))^2-(F''_{\dot{\eta}_k\dot{\xi}_k}(Z))^2\right|$$ $$=\left|\left(\frac{(\xi_k-u_k)^2-(\eta_k-v_k)^2}{((\xi_k-u_k)^2+(\eta_k-v_k)^2)^2}-2+2S_1-2S_3\right)^2\right.$$ $$\left.+\left(\frac{4(\xi_k-u_k)(\eta_k-v_k)}{((\xi_k-u_k)^2+(\eta_k-v_k)^2)^2}+2S_2-2S_4\right)^2\right|,$$
where
$$S_3=\sum_{i=1}^n\frac{(\xi_k-x_i)^2-\eta_k^2}{((\xi_k-x_i)^2+\eta_k^2)^2} \ws \ws \mbox{and} \ws \ws S_4=\sum_{i=1}^n\frac{2(\xi_k-x_i)\eta_k}{((\xi_k-x_i)^2+\eta_k^2)^2}.$$
Let us observe that
$$\frac{\left|(\xi_k-x_i)^2-\eta_k^2\right|}{((\xi_k-x_i)^2+\eta_k^2)^2}, \frac{2\left|(\xi_k-x_i)\eta_k\right|}{((\xi_k-x_i)^2+\eta_k^2)^2}\leq \frac{1}{(\xi_k-x_i)^2+\eta_k^2}.$$
To estimate $S_3$ and $S_4$ we show that for any $a$ and $b\neq 0$ real constants
\begin{equation}\sum_{i=1}^n\frac{1}{(a-x_i)^2+b^2}=O(\sqrt{n}).\end{equation}
Indeed, according to \cite[Corollary 4.3]{km}, if $c_1<, \dots <c_{m+n}$ are the zeros of the classical Hermite polynomial $h_{m+n}$, then at least $n-r$ intervals $(c_k,c_{k+1}), 1\leq k<n+m$ contain a (regular) zero of $P_{m+n}$, where $r$ is length of $\lambda$. Let these intervals be $I_1,\dots ,I_{n_1}$, $n-r\leq n_1\leq n$ listed in increasing order. Let us choose from every interval $I_k$ one zero of $P_{m+n}$, $x_{i_k}$. Let us denote by $\mathcal{I}$ the collection of the chosen indices and by $\mathcal{I}_e$ and $\mathcal{I}_o$ the chosen even ($i_2,i_4, \dots$) and odd ($i_1,i_3, \dots$) indices respectively. Since by \cite[6.31.22]{sz} $\frac{c}{n^{\frac{1}{6}}}>c_{k+1}-c_k>\frac{c}{\sqrt{n}}$ and $c_{i_k+r+m+1}-c_{i_k}>\Delta _{i_k}>c_{i_k+2}-c_{i_k+1}$ thence the previous estimation holds for $\Delta _{i_k}$ as well, where  $\Delta _{i_k}:=x_{i_{k+2}}-x_{i_k}$, and we have
$$\sum_{i=1}^n\frac{1}{(a-x_i)^2+b^2}=\sum_{i\notin \mathcal{I}}(\cdot)+\sum_{i_k\in \mathcal{I}_e}(\cdot)+\sum_{i_k\in \mathcal{I}_o}(\cdot)$$
$$\leq c(b,r)+c \sqrt{n} \left(\sum_{i_k\in \mathcal{I}_e}\frac{\Delta _{i_k}}{(a-x_{i_k})^2+b^2}+\sum_{i_k\in \mathcal{I}_o}\frac{\Delta _{i_k}}{(a-x_{i_k})^2+b^2}\right)=O(\sqrt{n}).$$
Let $\varrho:=\varrho(k,n)=\sqrt{(\xi_k-u_k)^2+(\eta_k-v_k)^2)}$. By \cite[Theorem 2.3.]{km} $\varrho\leq \frac{c}{\sqrt{n}}$. Thus according to (5)

$$|\det A_{k,k}|$$ $$\geq c \left|\left(\frac{(\xi_k-u_k)^2-(\eta_k-v_k)^2}{((\xi_k-u_k)^2+(\eta_k-v_k)^2)^2}\right)^2+\left(\frac{4(\xi_k-u_k)(\eta_k-v_k)}{((\xi_k-u_k)^2+(\eta_k-v_k)^2)^2}\right)^2-c\frac{\sqrt{n}}{\varrho^2}\right|$$ $$\geq  \frac{c_1}{\varrho^4}-c\frac{\sqrt{n}}{\varrho^2}>c\frac{n}{\varrho^2}.$$
Finally, by the previous calculations
$$\left|F''_{\dot{\xi}_k\dot{\xi}_k}(Z)\right|, \left|F''_{\dot{\eta}_k\dot{\xi}_k}(Z))\right|<cn,$$
which finishes the proof of (3).

The proof of (4) is similar, so we omit the details.

\medskip

The main tool of proving Theorems 1 and 2 is some lemmas on partitioned matrices. Following the definitions of R. Varga (cf. \cite{v}), by a partition $\pi$ of $\mathbb{C}^n$, we mean a finite collection $\{W_i\}_{i=1}^l$ of pairwise disjoint linear subspaces, each having dimension at least unity, whose direct sum is $\mathbb{C}^n$. Furthermore let $\Phi_{\pi}:=(\Phi_1, \dots , \Phi_l)$, where $\Phi_j$ is a norm on $W_j$ , for each $j=1, \dots , l$. $A = [A_{i,j} ] \in \mathbb{C}^{n\times n}$ is strictly block diagonally dominant with respect to $\Phi_{\pi}$ if
$$\left(\|A_{ii}^{-1}\|_{\Phi_{\pi}}\right)^{-1}>\sum_{1\leq j\leq l \atop j\neq i}\|A_{ij}\|_{\Phi_{\pi}}, \ws \ws 1\leq i\leq l.$$

\begin{lemma}\cite[Theorem 6.2.]{v}
Given a partition $\pi$ of $\mathbb{C}^n$ and given $\Phi_{\pi}$, assume that
$A = [A_{i,j} ] \in \mathbb{C}^{n\times n}$, partitioned by $\pi$, is strictly block diagonally dominant
with respect to $\Phi_{\pi}$. Then, $A$ is nonsingular.\end{lemma}

\medskip

\noindent {\it Proof of Theorem 1.}
We prove the second statement by the lemma above. The partition is as follows. We have $m$ groups of variables with two members,\\ $(\xi_1,\eta_1),\dots , (\xi_m,\eta_m)$, and $n$ groups with a single member, $x_1, \dots , x_n$. That is for $k,l\in \{1,\dots , m\}$, $i,j\in \{1,\dots , n\}$
$$A_{k,l}=\left[\begin{array}{cc}h_{2k-1,2l-1}&h_{2k-1,2l}\\h_{2k,2l-1}&h_{2k,2l}\end{array}\right]; \ws \ws A_{k,m+i}=\left[\begin{array}{c}h_{2k-1,2m+i}\\h_{2k,2m+i}\end{array}\right]; $$ $$\ws \ws A_{m+i,l}=\left[\begin{array}{cc}h_{2m+i,2l-1}&h_{2m+i,2l}\end{array}\right]; \ws \ws A_{m+i,m+j}=h_{2m+i,2m+j}.$$
Lemma 1 ensures that the first $2m$ rows ($m$ groups) of {\bf H} are strongly block diagonally dominant, that is
$$\|A_{k,k}^{-1}\|^{-1}> c n>>c\sqrt{n} >\sum_{1\leq l\leq m \atop l\neq k} \|A\|_{k,l}+\sum_{i=1}^n\|A_{km+i}\|, \ws\ws k=1, \dots , m.$$
The remainding $n$ rows are not necessarily block diagonally dominant, but we know that $\sum_{l=1}^m\|A_{m+i,l}\|=c(m)$ that is the sum depends on $m$, $2\sum_{k=1}^m\frac{(x_i-u_k)^2-v_k^2}{((x_i-u_k)^2+v_k^2)^2}-\frac{(\xi_k-x_i)^2-\eta_k^2}{((\xi_k-x_i)^2+\eta_k^2)^2}$ tends to zero, when $n$ tends to infinity (cf. \cite{km}). We multiply $A_{m+i,l}$ ($i=1,\dots n; l=1, \dots, m$) by a constant $\frac{1}{K(m)}$ which depends on $m$ but which is independent of $n$, and parallely multiply $A_{k,m+i}$ ($i=1,\dots n; k=1, \dots, m$) by $K(m)$, that is we transform ${\bf H}$ by a diagonal matrix $D$ with $d_{kk}=1$, $k=1, \dots, 2m$, $d_{2m+i2m+i}=K(m)$, $i=1,\dots n$.
Since $n$ is large enough, it can be archived that the new matrix ${\bf \hat{H}}=D^{-1}{\bf H}D$ is strictly block diagonally dominant, and its eigenvalues necessarily coincide with the eigenvalues of ${\bf H}$. Thus by Lemma 2 ${\bf H}$ is nonsingular.

\medskip

Now we turn to the proof of Theorem 2. To this end we have to cite another theorem of R. Varga on the location of the eigenvalues of partitioned matrices. The Gersgorin set with respect to the partition $\pi$ and the norm $\Phi_{\pi}$ is $G_{\pi}^{\Phi_{\pi}}(A)=\cup_{i=1}^lG_{i,\pi}^{\Phi_{\pi}}(A)$, where $G_{i,\pi}^{\Phi_{\pi}}(A):=\left\{z\in \mathbb{C} : \left(\|(zI_i-A_{i,i})^{-1}\|_{\Phi_{\pi}}\right)^{-1}\leq \sum_{1\leq j\leq l \atop j\neq i}\|A_{ij}\|_{\Phi_{\pi}}\right\}$, ($1\leq i\leq l)$, where $I_i$ denotes the identity matrix for the subspace $W_i$. Let $\sigma(A)$ be the spectrum of $A$.

\medskip

\begin{lemma}\cite[Theorem 6.3.]{v}
Given a partition $\pi$ of $\mathbb{C}^n$ and given $\Phi_{\pi}$, let $A = [A_{i,j} ] \in \mathbb{C}^{n×n}$, partitioned by $\pi$. If $\lambda \in \sigma(A)$, there is an $i \in \{1, \dots , l\}$ such that $\lambda \in G_{i,\pi}^{\Phi_{\pi}}(A)$. That is $\sigma(A)\subset G_{\pi}^{\Phi_{\pi}}(A)$.\end{lemma}

\medskip

\noindent {\it Proof of Theorem 2.}
With the partition and norm above, let
\begin{equation} G_r=\left(\cup_{i=1}^nG_{i,\pi}^{\Phi_{\pi}}({\bf \hat{H}})\right)\cap \mathbb{R},\end{equation}
the union of the Gersgorin discs around $h_{2m+i,2m+i}$, and
\begin{equation} G_e=\left(\cup_{k=1}^mG_{k,\pi}^{\Phi_{\pi}}({\bf \hat{H}})\right)\cap \mathbb{R},\end{equation}
the union of the Gersgorin sets with respect to $A_{k,k}$.
Since ${\bf \hat{H}}$ is strictly block diagonally dominant, and the "blocks" $A_{m+i,m+i}\in\mathbb{R}_-$, the first part of the statement on $G_r$ is immediately proved. Recalling the computations in the proof of Theorem 1, the estimation on modulus of the eigenvalues in $G_r$ is also obvious.

According to the previous lemma,
$$G_e=\cup_{k=1}^m G_{e,k}=\cup_{k=1}^m\left\{x\in\mathbb{R} : \left(\|xI_2-A_{k,k}\|^{-1}\right)^{-1}\leq R\right\},$$
where $R=\sum_{1\leq l\leq m \atop l\neq k}\|A_{k,l}\|$. Since $F''_{\dot{\eta}_k\dot{\eta}_k}(Z)=-F''_{\dot{\xi}_k\dot{\xi}_k}(Z)$,
$$\left(\|xI_2-A_{k,k}\|^{-1}\right)^{-1}=\frac{\left|x^2-\left(\left(F''_{\dot{\xi}_k\dot{\xi}_k}\right)^2(Z)+\left(F''_{\dot{\eta}_k\dot{\xi}_k}\right)^2(Z)\right)\right|}{|x|+|F''_{\dot{\xi}_k\dot{\xi}_k}(Z)|+|F''_{\dot{\eta}_k\dot{\xi}_k}(Z)|}=:\frac{\left|x^2-U^2\right|}{|x|+V}.$$
Thus if $|x|\geq U$,
$$x\in G_{e,k} \ws \ws \mbox{iff } \ws \ws \left(|x|-\frac{R}{2}\right)^2\leq RV+U^2+\frac{R^2}{4},$$
and if $|x|\leq U$,
$$x\in G_{e,k} \ws \ws \mbox{iff } \ws \ws \left(|x|+\frac{R}{2}\right)^2\geq U^2-RV+\frac{R^2}{4}.$$
So if $x\in G_{e,k}$, then
$$U\leq |x|\leq \sqrt{RV+U^2+\frac{R^2}{4}}+\frac{R}{2}, \ws\ws \mbox{or} \ws\ws \sqrt{U^2-RV+\frac{R^2}{4}}-\frac{R}{2}\leq |x|\leq U.$$
It is easy to see that the two sets above are not empty. Recalling that  $U^2> cn^2$, $R=O(\sqrt{n})$ and $V=O(n)$ the remainding part of the proof is obvious.

\medskip

For the proof of Theorem 3 we need some lemmas. First we give the differential equation of $d_{\nu}$.

\medskip

\lemma
For $\nu=0,1,\dots$, $d_{\nu}$ fulfils the following differential equation
$$d_{\nu}''-\left(2x+2\frac{h'_{\nu}}{h_{\nu}}\right)d_{\nu}'+4x\frac{h'_{\nu}}{h_{\nu}}d_{\nu}=0.$$

\medskip

\proof
Let us consider the Schr\"odinger operator $-\varphi''+V(x)\varphi$ with eigenfunctions and eigenvalues $\varphi_j=h_je^{-\frac{x^2}{2}}$ and $E_j$ respectively.
Let
\begin{equation}\bar{\varphi}_m=\frac{\mathrm{Wr}[\varphi_k,\varphi_m]}{\varphi_k}=\frac{e^{-\frac{x^2}{2}}}{h_k}d_{k,m},\end{equation}
where
$$d_{k,m}=\left|\begin{array}{cc}h_{k}&h_{m}\\h'_{k}&h'_{m}\end{array}\right|.$$
According to \cite[(27)]{gfgu}
\begin{equation}-\bar{\varphi}_m''+(V-\left(\log\varphi_k)''\right)\bar{\varphi}_m=E_m\bar{\varphi}_m.\end{equation}
Substituting the right-hand side of (8) to (9) we have
$$d_{k,m}''-(\left(2x+2\frac{h'_{k}}{h_{k}}\right)d_{k,m}'+2\left(2x\frac{h_k'}{h_{k}}+m-k-1\right)d_{k,m}=0,$$
which gives the statement.

\medskip

Now we can give the elements of the Hessian with respect to $d_{\nu}$ in a simple form.

\medskip

\lemma
$$h_{2k-1,2k-1}=F''_{\dot{\xi}_k\dot{\xi}_k}(Z)=-h_{2k,2k}=F''_{\dot{\eta}_k\dot{\eta}_k}(Z)=\Re r_{m,n}(z_k),$$ $$ h_{2k-1,2k}=F''_{\dot{\eta}_k\dot{\xi}_k}(Z)=\Im r_{m,n}(z_k),$$
where
$$r_{m,n}(x)=-\left(8\left(x^2+1-\sqrt{\frac{\nu+1}{2}}\frac{2h_{\nu}^2+h_{\nu+1}^2}{d_{\nu}}\right)+2n\right).$$

\medskip

\proof
\begin{equation}-h_{2k-1,2k-1}=\Re\left(2+2\left(\frac{d_{\nu}'}{d_{\nu}}\right)'(z_k)-\left(\frac{P_{m,m+n}'}{P_{m,m+n}}\right)'(z_k)\right).\end{equation}
By the differential equations of $P_{m,m+n}$ and $d_{\nu}$ we have
$$\left(\frac{P_{m,m+n}'}{P_{m,m+n}}\right)'(z_k)=\left(2x+2\frac{d_{\nu}'}{d_{\nu}}-\frac{P_{m,m+n}}{P_{m,m+n}'}\left(\left(4x+2\frac{h_{\nu}'}{h_{\nu}}\right)\frac{d_{\nu}'}{d_{\nu}}-4x\frac{h_{\nu}'}{h_{\nu}}+2n\right)\right)'|_{x=z_k}$$
$$=2+2\left(\frac{d_{\nu}'}{d_{\nu}}\right)'(z_k)-\left(\left(4x+2\frac{h_{\nu}'}{h_{\nu}}\right)\frac{d_{\nu}'}{d_{\nu}}-4x\frac{h_{\nu}'}{h_{\nu}}+2n\right)(z_k).$$
Thus
$$-r_{m,n}(z_k)=\left(\left(4x+2\frac{h_{\nu}'}{h_{\nu}}\right)\frac{d_{\nu}'}{d_{\nu}}-4x\frac{h_{\nu}'}{h_{\nu}}+2n\right)(z_k).$$
Taking into consideration that
$$\frac{d_{\nu}'}{d_{\nu}}=2x-2\frac{h_{\nu}h_{\nu+1}}{d_{\nu}} \ws \mbox{and} \ws
\frac{h_{\nu}'}{h_{\nu}}\frac{d_{\nu}'}{d_{\nu}}=2x\frac{h_{\nu}'}{h_{\nu}}-2\frac{h_{\nu}'h_{\nu+1}}{d_{\nu}}$$
we have
$$-r_{m,n}(z_k)=8x^2+2n-\frac{8xh_{\nu}h_{\nu+1}+4h_{\nu}'h_{\nu+1}}{d_{\nu}}.$$
Finally observing that
$$2xh_{\nu}h_{\nu+1}=\sqrt{2(\nu+1)}\left(h_{\nu}^2+h_{\nu+1}^2\right)-d_{\nu}\ws \mbox{and} \ws h_{\nu}'h_{\nu+1}=\sqrt{2(\nu+1)}h_{\nu}^2-d_{\nu}$$
the proof is finished.

\medskip

As it is mentioned above $d_{\nu}$ has simple roots, that is by \cite[Theorem 2.3]{km} $|z_k-w_k|\leq \frac{c}{\sqrt{n}}$, $k=1,\dots ,m$. We show that the distance of the zeros of $d_{\nu}$ and of $P_{m,m+n}$ cannot be too small.

\medskip

\lemma  Let $H$, $P_{m,m+n}$ be the same as above, $w_k$, $z_k$ ($k=1, \dots , m$) the zeros of $H$ and the exceptional zeros of $P_{m,m+n}$ respectively. If $n$ is large enough, there is a $c>0$ such that for $k=1,\dots ,m$ $|z_k-w_k|\geq \frac{c}{\sqrt{n}\log n}$.

\medskip

\proof
By \cite[(5.9)]{km}
\begin{equation}\frac{1}{w_k-z_k}=w_k+\sum_{1\leq l\leq m \atop l\neq k}\frac{1}{w_k-w_l}-\sum_{1\leq l\leq m \atop l\neq k}\frac{1}{w_k-z_l}-\sum_{j=1}^n\frac{1}{w_k-x_j}.\end{equation}
The first two terms of the right-hand side obviously depend only on $m$, and the third term by \cite[Theorem 2.3]{km} can be estimated by a constant depending only on $m$. We estimate the modulus of the $4^{th}$ term. As in the proof of \cite[Theorem 2.3]{km} we have
$$\left|\Im \sum_{j=1}^n\frac{1}{w_k-x_j}\right|=|v_k|\sum_{j=1}^n\frac{1}{|w_k-x_j|^2}=O(\sqrt{n}),$$
cf. (5). Similarly to (5)
$$\left|\Re \sum_{j=1}^n\frac{1}{w_k-x_j}\right|\leq \frac{1}{v_k^2}\sum_{j\notin \mathcal{I} \ws \mathrm{or} \atop j, |u_k-x_j|\leq 1}1+c \sqrt{n}\left| \sum_{j\in \mathcal{I}_e \atop |u_k-x_j|>1}\frac{(u_k-x_j)\Delta_j}{(u_k-x_j)^2+v_k^2}\right|$$ $$+
c \sqrt{n}\left| \sum_{j\in \mathcal{I}_o \atop |u_k-x_j|>1}\frac{(u_k-x_j)\Delta_j}{(u_k-x_j)^2+v_k^2}\right|\leq c\left(\sqrt{n}+\sqrt{n}\log n+\sqrt{n}\log n\right),$$
where we used that the number of zeros on $[a,b]$ is $O(\sqrt{n})$ cf. \cite{km}, and the estimation  $|u_k-x_j|<c\sqrt{n}$ ($j\in\mathcal{I}$), cf. \cite[6.32.6]{sz}.

\medskip

\noindent{\it Proof of Theorem 3.} Let us recall that
$$-r_{m,n}(z_k)=2n+8(z_k^2+1)-4\sqrt{2(\nu+1)}\frac{\left(2h_{\nu}^2+h_{\nu+1}^2\right)(z_k)}{\prod_{1\leq l \leq m \atop l\neq k}(z_k-w_l)}\frac{1}{z_k-w_k}.$$
Let us observe that the second term depends on $m$ and the same is valid for the first factor of the third term. Thus the modulus of the third term is $O(\sqrt{n}\log n)$. That is
$$|h_{2k-1,2k-1}|=|h_{2k,2k}|\sim n \ws \mbox{and} \ws |h_{2k-1,2k}|=o(n).$$
Comparing this with the proof of Theorem 1 the proof is finished.

\medskip

\noindent{\it Proof of Proposition 1.}
Obviously, as in \cite{ahor}
$$\frac{\partial \log T_w}{\partial u_i}(Z)=0, \ws\ws i=1, \dots , n.$$
With the notations above the off-diagonal elements of the Hessian are
$$h_{i,j}=\frac{2}{(x_i-x_j)^2}$$
and
$$h_{i,i}=(\log w)''(x_i)-2\sum_{1\leq i<j\leq n}\frac{1}{(x_i-x_j)^2}.$$
Recalling \cite[Proposition 2]{ho} for uniqueness, the proof is complete.

\medskip

\noindent{\it Proof of Proposition 2.}
Since $P_{m,m+n}=P_mq_n$ fulfils the differential equation (1):
$$y''+M_ny'+N_ny=0,$$
$q_n$ fulfils the following differential equation:
$$y''+M_{1,n}y'+N_{1,n}y=0,$$
where
$$M_{1,n}=M_n+2\frac{P_m'}{P_m}, \ws\ws N_{1,n}=N_n+M_n\frac{P_m'}{P_m}+\frac{P_m''}{P_m},$$
(cf. \cite{ahor})
Since $M_{1,n}=\left(\log w_1\right)'$ it is enough to show that  $M_{1,n}'<0$.
$$M_{1,n}'=\left(M_n+2\frac{P_m'}{P_m}\right)'=M_{0,n}'+2\left(\frac{P_m'}{P_m}-\frac{H'}{H}\right)'$$ $$=M_{0,n}'-2\sum_{k=1}^m\left(\frac{1}{(x-z_k)^2}-\frac{1}{(x-w_k)^2}\right).$$
Thus with the notations above, if $z_k \to w_k$ $k=1,\dots ,m$ when $n \to \infty$, if $n$ is large enough $w_1$ is admissible if $w_0$ is admissible. Applying the previous lemma the proof is finished.

\medskip

\noindent{\it Proof of Corollary.}
Taking into consideration that $M_{0,n}'=-2$, by \cite[Theorem 2.3]{km} and by Proposition 2 the proof is complete.

\medskip

We remark here that the differential equation of $q_n$ has the following form:
$$y''+\left(-2x-2\frac{H'}{H}+2\frac{P_m'}{P_m}\right)y'$$ $$+\left(\frac{P_m''}{P_m}+\frac{P_m'}{P_m}\left(-2x-2\frac{H'}{H}\right)+\frac{H''}{H}+2x\frac{H'}{H}+2n-m\right)y=0.$$

\medskip

\noindent \small{Department of Analysis, \newline
Budapest University of Technology and Economics}\newline
\small{ g.horvath.agota@renyi.mta.hu}

\end{document}